\documentclass[11pt,a4paper,reqno]{amsart}
\usepackage[T1]{fontenc}
\usepackage[utf8]{inputenc}
\usepackage[nomath]{lmodern}
\usepackage[english]{babel}
\usepackage{amsmath,amssymb,amsthm,mathtools}
\usepackage[margin=25mm]{geometry}
\usepackage{microtype,needspace}
\usepackage[unicode,colorlinks=true,allcolors=blue]{hyperref}
\hypersetup{
  pdftitle={Division ring extensions and the realization of Coxeter diagrams I2(n)},
  pdfauthor={A. I. Generalov, A. Tursunbaev}
}
\numberwithin{equation}{section}
\newtheorem{theorem}{Theorem}[section]
\newtheorem{lemma}[theorem]{Lemma}

\theoremstyle{remark}
\newtheorem{remark}[theorem]{Remark}
\newcommand{\Frac}{\operatorname{Frac}}
\newcommand{\Div}{\operatorname{Div}_k}

\newcommand{\Hom}{\operatorname{Hom}}
\newcommand{\Q}{\mathbb Q}
\newcommand{\Z}{\mathbb Z}
\makeatletter
\renewcommand{\@@and}{\unskip,}
\def\@settitle{\begin{center}\normalfont\large\bfseries\@title\end{center}}
\renewcommand{\@setkeywords}{%
  {\itshape\keywordsname:}\enspace\@keywords\@addpunct.}
\makeatother
\addto\captionsenglish{\renewcommand{\keywordsname}{Keywords}}
\title[Division ring extensions and Coxeter diagrams]
{Division ring extensions and the realization\\
of Coxeter diagrams \texorpdfstring{$I_2(n)$}{I2(n)}}

\author{A.~I.~Generalov}
\address{St.~Petersburg State University,
Universitetskii pr.~28,
St.~Petersburg 198504, Russia}
\email{\href{mailto:ageneralov@gmail.com}{ageneralov@gmail.com}}

\author{A.~Tursunbaev}
\address{St.~Petersburg State University,
Universitetskii pr.~28,
St.~Petersburg 198504, Russia}
\email{\href{mailto:tursunbaev.art@gmail.com}{tursunbaev.art@gmail.com}}

\date{}

\begin{document}

\begin{abstract}
For each integer $n\ge5$, we construct a hereditary artinian ring
of finite representation type with Coxeter diagram $I_2(n)$. Both of its diagonal division rings are two-generated over the field $\Q(\zeta)$,
where $\zeta$ is a primitive root of unity of order $n-2$.
The construction is based on a chain of quantized Weyl relations
and a shift automorphism. We give one-sided bases
for the corresponding division ring extensions and compute the right
dimension sequence of the bimodule. The resulting ring has exactly $n$ isomorphism classes of indecomposable finitely generated right modules.
\end{abstract}
\maketitle

\section{Introduction}
\label{sec:introduction}

Dowbor, Ringel, and Simson~\cite{DRS} associated Coxeter diagrams with
hereditary artinian rings of finite representation type.
The realization of such diagrams is related to the construction of bimodules
over division rings subject to conditions on the dimensions of their
successive dual bimodules. Schofield~\cite{Sch}
constructed an example with diagram $I_2(5)$.
A criterion for realizing dimension sequences in terms of the dimensions
of bimodules defined by stabilizers of subspaces in a division ring
extension was given by \v{S}aroch~\cite[Lemmas~1.3 and~1.4]{Sar}.

In this paper, we give an explicit construction of algebras with diagrams
$I_2(n)$ for all integers $n\ge5$ in characteristic zero.
The construction is parametrized by a primitive root of unity of order $N=n-2$.
We use quantized Weyl relations on a doubly infinite chain
and the automorphism that shifts its indices.
The relations on finite intervals of the chain are related to the linear
connected quantized Weyl algebras studied by Fish and Jordan~\cite{FishJordan}.
Passing to a skew field of rational functions yields two-generated
diagonal division rings and allows us to express explicitly the generators
of one in terms of those of the other.

Our main result is as follows.

\begin{theorem}\label{thm:quantum}
Let $n\ge5$, $N=n-2$, and let $\zeta$ be a primitive $N$th root of unity.
Over the field $k=\Q(\zeta)$, there exist two-generated division rings
$H\subset E$ such that
\[
 \dim(E_H)=2,\qquad
 \dim({}_H E)=N,
\]
and the right dimension sequence of the bimodule
${}_H E_E$ is periodic with period
\begin{equation}\label{eq:qperiod}
 \left(1,\underbrace{2,\ldots,2}_{N-1},1,N\right).
\end{equation}
The algebra
\[
 \Lambda_n=\begin{pmatrix}H&E\\0&E\end{pmatrix}
\]
is hereditary, left and right artinian, and has exactly $n$
isomorphism classes of indecomposable finitely generated right modules.
Its Coxeter diagram is $I_2(n)$.
\end{theorem}

All division rings in this paper are associative and are not assumed to be commutative.
For a subset $U$ of a division ring containing the base field $k$,
let $\Div(U)$ denote the smallest division subring containing
$k$ and $U$. Thus generation is understood in the sense of division rings,
allowing inverses.
For an $(F,G)$-bimodule $M$, its right dual bimodule is
$M^R=\Hom_G(M_G,G_G)$ with its natural
$(G,F)$-bimodule structure. The right dimension sequence of $M$
consists of the right dimensions of $M,M^R,M^{RR},\ldots$
over the corresponding division rings.

In Section~\ref{sec:chain}, we construct the division rings $E$ and $G_r$.
In Section~\ref{sec:bases}, we establish the right linear independence
of powers of $t$ and the decompositions
\[
 E=\bigoplus_{i=0}^{r}t^iG_r
        =\bigoplus_{i=0}^{N-r}G_rt^i
 \qquad(1\le r<N).
\]
In Section~\ref{sec:dimensions}, we compute the stabilizers
of the corresponding subspaces and complete the proof
of Theorem~\ref{thm:quantum}.

\section{The quantum chain and division ring extensions}
\label{sec:chain}

Fix $N\ge3$ and let $\zeta$ be a primitive $N$th root of unity; put $k=\Q(\zeta)$.
For a natural number $a$, put
\[
 [a]_\zeta=\frac{\zeta^a-1}{\zeta-1}.
\]
Then $[a]_\zeta \ne0$ for $1\le a<N$ and $[N]_\zeta =0$.
Consider the $k$-algebra $A $ with generators $x_i$, $i\in\Z$,
and relations
\begin{equation}\label{eq:qchain}
 x_jx_i=\zeta^{(-1)^{j-i+1}}x_ix_j+\delta_{j,i+1}
 \qquad (j>i),
\end{equation}
where $\delta$ is the Kronecker delta. On a finite interval, our relations are obtained from those of the linear connected quantized Weyl algebra in~\cite[Example~2.6]{FishJordan} by taking $q=\zeta^{-1}$ and multiplying the generators in positions $2,4,6,\ldots$ by $(1-\zeta)^{-1}$.

\begin{lemma}\label{lem:qore}
The algebra $A $ is a left and right Ore domain.
Its ordered monomials form a $k$-basis, and the map
$x_i\mapsto x_{i+1}$ extends to an automorphism $\beta$ of its
division ring of fractions $T =\Frac(A )$.
\end{lemma}

\begin{proof}
For a finite interval $[a,b]$, successively adjoin
$x_a,\ldots,x_b$. When adjoining $x_j$, use the Ore extension
with automorphism and skew derivation
\[
 \sigma_j(x_i)=\zeta^{(-1)^{j-i+1}}x_i,\qquad
 \delta_j(x_i)=
 \begin{cases}1,&i=j-1,\\0,&i<j-1.\end{cases}
\]
Here $\delta_j(uv)=\sigma_j(u)\delta_j(v)+\delta_j(u)v$.
The map $\sigma_j$ preserves all the preceding relations:
for an adjacent pair, the product of the two scalar factors is~$1$.
For $\delta_j$, it suffices to check the relations involving $x_{j-1}$.
If $i<j-1$, then, for $c=\zeta^{(-1)^{j-i}}$, we have
$\delta_j(x_{j-1}x_i-cx_ix_{j-1})=x_i-c\sigma_j(x_i)=0$.
We therefore obtain an iterated Ore extension, which is a noetherian
domain with the stated basis.

The embeddings associated with finite intervals are compatible with these bases.
Any two elements of $A $ lie in the algebra on a single finite interval,
where both Ore conditions hold. Hence both conditions also hold in $A $.
The shift of indices preserves~\eqref{eq:qchain} and is invertible,
so it extends to an automorphism of $T $.
\end{proof}

Put
\[
 E =T (s;\beta),\qquad sa=\beta(a)s\quad(a\in T ).
\]
Here $T (s;\beta)$ is the division ring of fractions of the skew polynomial ring $T [s;\beta]$.
Introduce the notation
\begin{equation}\label{eq:qgens}
 \begin{gathered}
 x=x_{-1},\qquad t=x_0,\qquad y=x_1,\\
 v_r=y^rs=st^r\ (r\ge0),\qquad
 G_r=\Div(x,v_r)\ (1\le r<N).
 \end{gathered}
\end{equation}
The identities $x_i=s^{i+1}xs^{-i-1}$ imply $E =\Div(x,s)$.
Relations~\eqref{eq:qchain} give
\begin{equation}\label{eq:qthree}
 tx=\zeta xt+1,\qquad yt=\zeta ty+1,\qquad yx=\zeta^{-1}xy.
\end{equation}
For $r\ge1$, these yield
\begin{equation}\label{eq:qvr}
 v_rt=v_{r+1},\qquad
 v_rx-\zeta^rtv_r=[r]_\zeta v_{r-1},
 \qquad v_0=s.
\end{equation}
Moreover, each element $x_i^N$ is central in $T $.
Indeed, for an adjacent pair this follows from the formula
$ba^h=\zeta^ha^hb+[h]_\zeta a^{h-1}$ when $ba=\zeta ab+1$,
and for all other pairs it follows directly from~\eqref{eq:qchain}.
In particular,
\begin{equation}\label{eq:qpower}
 v_rx^Nv_r^{-1}=t^N\in G_r\qquad(1\le r<N).
\end{equation}

\section{Left and right bases}
\label{sec:bases}

\begin{lemma}\label{lem:qindependence}
If $1\le r<N$, then the elements
$1,t,\ldots,t^r$ are right linearly independent over $G_r$.
\end{lemma}

\begin{proof}
Put $v=v_r$, $P_0=1$, and $P_b=x_1^r\cdots x_b^r$ for $b\ge1$.
Let $B_j$ denote the division subring of $T $ generated over $k$
by the elements $x_i$ with $i\le j$. Let $\theta$ be conjugation
by $v$ on $T $, and put $K=\Div(\theta^a(x)\mid a\in\Z)$.
Then
\begin{equation}\label{eq:qunion}
 K=\bigcup_{b\ge1}P_bB_{b-1}P_b^{-1}.
\end{equation}
Indeed, for $a\le0$, the element $\theta^a(x)$ differs
from $x_{a-1}$ by a nonzero factor in $k$.
For $a\ge1$, the equality $v^a=P_as^a$ gives
$\theta^a(x)=P_ax_{a-1}P_a^{-1}$.
If $1\le a\le b$, then replacing $P_a$ by $P_b$ in the latter
expression changes it only by a nonzero scalar factor: the indices of
the additional generators differ from the index of $x_{a-1}$ by at least two.
For the same reason, conjugation by $P_b$ multiplies
each $x_i$, $i<0$, by a nonzero scalar. Hence
\[
 P_bB_{b-1}P_b^{-1}=\Div(\theta^a(x)\mid a\le b),
\]
which proves~\eqref{eq:qunion} and shows that the division subrings in this union form an increasing chain.

By the definition of $K$, we have $\theta(K)=K$ and $K\subseteq G_r$,
since $\theta^a(x)=v^axv^{-a}$ for all $a\in\Z$.
The equalities $v^j=P_js^j$ and the left linear independence of the powers of $s$
over $T$ imply that the powers of $v$ are left linearly independent over $T$. Thus the subring of $E$ generated by $K$ and $v$ is a skew polynomial ring with the relation $va=\theta(a)v$ for $a\in K$. Since the restriction of $\theta$ to $K$ is an automorphism of the division ring $K$, this ring satisfies the right Ore condition. Its right fractions form a division subring of $E$ equal to $G_r$: it contains $x,v$ and, on the other hand, is contained in $G_r$ because $K\subseteq G_r$.

We show that it suffices to prove the right linear independence of
$1,t,\ldots,t^r$ over $K$. Suppose that
\[
 \sum_{i=0}^r t^ia_i=0,\qquad a_i\in G_r.
\]
Choose a common right denominator $q$ in the above polynomial ring and write
\[
 a_iq=\sum_{j=0}^{m}b_{ij}v^j,\qquad b_{ij}\in K
 \quad(0\le i\le r).
\]
Multiplying the original equality on the right by $q$, we obtain
\[
 \sum_{j=0}^{m}\left(\sum_{i=0}^r t^ib_{ij}\right)v^j=0.
\]
The coefficients in parentheses belong to $T$. Hence the left linear
independence of the powers of $v$ over $T$ gives
$\sum_{i=0}^r t^ib_{ij}=0$ for every $j$.
If $1,t,\ldots,t^r$ are right linearly independent over $K$,
then all $b_{ij}=0$, whence $a_iq=0$ and $a_i=0$.

In a relation over $K$, all coefficients belong, by~\eqref{eq:qunion},
to a single division subring $P_bB_{b-1}P_b^{-1}$, since
these division subrings form an increasing chain. Writing the coefficients
as $P_bc_iP_b^{-1}$ and multiplying the relation on the right by $P_b$,
we reduce the proof to the implication
\begin{equation}\label{eq:qdependence}
 \sum_{i=0}^r t^iP_b c_i=0,\quad c_i\in B_{b-1}
 \quad\Longrightarrow\quad c_0=\cdots=c_r=0.
\end{equation}

On $B_b$ we use the skew derivation $\delta$
from the extension obtained by adjoining $x_{b+1}$.
It extends to the division ring of fractions, vanishes on $B_{b-1}$,
and the corresponding automorphism satisfies
$\sigma(t)=\zeta^{(-1)^b}t$, $\sigma(x_b)=\zeta x_b$.
Applying $\delta^\nu$ to~\eqref{eq:qdependence} for $0\le \nu\le r$
and dividing the rows by nonzero scalars, we obtain the system
$M_0c=0$. For $0\le m<b$, its subsequent transformations
are described by the matrices
\begin{equation}\label{eq:qmatrices}
 M_m[\nu,i]=\zeta^{\nu\varepsilon_m i}
 x_m^i x_{m+1}^r\cdots x_{b-1}^r x_b^{r-\nu},
 \quad \varepsilon_m=(-1)^{b-m},\quad 0\le \nu,i\le r.
\end{equation}
The division used to obtain $M_0$ is valid because
$[r]_\zeta [r-1]_\zeta \cdots[r-\nu+1]_\zeta \ne0$.

We prove that these matrices are invertible. First, let $m<b-1$,
$w=x_m$, $a=x_{m+1}$, and
$Y_\nu=a^rZ_\nu=x_{m+1}^r\cdots x_{b-1}^r x_b^{r-\nu}$.
Put $\chi=1$ if $b-m$ is odd, and $\chi=0$ otherwise.
Define commuting operators on column vectors by
\[
 (LU)_\nu=\zeta^{\nu\varepsilon_m}wU_\nu,\qquad
 (RU)_\nu=\zeta^{-r\chi}U_\nu w.
\]
Column $i$ of $M_m$ is $L^iY$.
Replace it by
\begin{equation}\label{eq:qcolumn}
 \prod_{h=0}^{i-1}(L-\zeta^hR)Y.
\end{equation}
Expanding the product shows that this is the original column
plus a right linear combination of the preceding columns.
Consequently, the simultaneous replacements~\eqref{eq:qcolumn}
are effected by right multiplication by an upper triangular matrix
with unit diagonal.

It follows from~\eqref{eq:qchain} that
\[
 wa^u=\zeta^{-u}a^uw-\zeta^{-u}[u]_\zeta a^{u-1},
 \qquad
 Z_\nu w=\zeta^{r(\chi-1)+\nu\varepsilon_m}wZ_\nu.
\]
Thus the $(\nu,i)$-entry of the transformed matrix equals
\begin{equation}\label{eq:qreduction}
 (-1)^i \zeta^{\nu i\varepsilon_m+i(i-1)/2-ri}
 \left(\prod_{h=0}^{i-1}[r-h]_\zeta \right)a^{r-i}Z_\nu.
\end{equation}
All scalar factors depending only on $i$ are nonzero.
Divide the columns by these factors, reverse the order of the columns, and multiply
row $\nu$ on the left by $\zeta^{-r\nu\varepsilon_m}$.
This yields $M_{m+1}$, since $\varepsilon_{m+1}=-\varepsilon_m$.
Thus $M_m$ is invertible if and only if $M_{m+1}$ is invertible.

For the last matrix, we have
\[
 M_{b-1}[\nu,i]=\zeta^{-\nu i}x_{b-1}^i x_b^{r-\nu}.
\]
Perform the same column replacements with $w=x_{b-1}$, $\chi=1$. The entry in position $(\nu,i)$ of the resulting matrix is
\[
 (-1)^i \zeta^{i(i-1)/2-ri}
 \left(\prod_{h=0}^{i-1}[r-\nu-h]_\zeta \right)x_b^{r-\nu-i}.
\]
Multiplying row $\nu$ on the left by $x_b^{\nu-r}$ and column $i$ on the right by
$(-1)^i\zeta^{ri-i(i-1)/2}x_b^i$, we obtain a matrix over $k$
\[
 C[\nu,i]=\prod_{h=0}^{i-1}[r-\nu-h]_\zeta,
 \qquad 0\le \nu,i\le r.
\]
If $i>r-\nu$, then the product contains the factor $[0]_\zeta=0$.
Moreover,
\[
 C[\nu,r-\nu]=\prod_{j=1}^{r-\nu}[j]_\zeta\ne0
 \qquad(0\le \nu\le r),
\]
where the empty product equals $1$. Hence, after reversing the order
of the columns, the matrix $C$ becomes upper triangular with nonzero
diagonal entries. Therefore $C$, and hence $M_0$,
is invertible.
This proves~\eqref{eq:qdependence}.
\end{proof}

\begin{lemma}\label{lem:qrightbasis}
For $1\le r<N$, we have the decomposition
\begin{equation}\label{eq:qrightbasis}
 E =\bigoplus_{i=0}^r t^iG_r.
\end{equation}
\end{lemma}

\begin{proof}
Let $U=\sum_{i=0}^rt^iG_r$.
Let $L_t$ denote left multiplication by $t$, and let $R_x$ denote
right multiplication by $x$. These operators commute.
From~\eqref{eq:qvr}, for $0\le j\le r$, we obtain
\[
 v_j=
 \frac{\prod_{h=j+1}^r(R_x-\zeta^hL_t)(v_r)}
 {[r]_\zeta [r-1]_\zeta \cdots[j+1]_\zeta }\in U.
\]
Here the products are empty when $j=r$. Since $st^j=v_j$,
we have $sU\subseteq U$. The formula
\[
 xt^j=\zeta^{-j}t^jx-\zeta^{-j}[j]_\zeta t^{j-1}
\]
gives $xU\subseteq U$.
The left stabilizer of the finite-dimensional right $G_r$-space $U$
is a division subring of $E $: the left multiplication on $U$ by a nonzero element of this stabilizer is injective and hence surjective. This stabilizer contains
$k,x,s$ and therefore equals $E $. Since $1\in U$, we obtain $U=E $.
The sum is direct by Lemma~\ref{lem:qindependence}.
\end{proof}

\begin{lemma}\label{lem:qadjacent}
If $1\le r<N-1$, then
\begin{equation}\label{eq:qadjacent}
 G_r+G_rt=G_{r+1}+tG_{r+1}.
\end{equation}
More generally, for $r\ge1$, $h\ge0$, and $r+h<N$, we have
\begin{equation}\label{eq:qflag}
 \sum_{i=0}^hG_rt^i=\bigoplus_{i=0}^h t^iG_{r+h}.
\end{equation}
\end{lemma}

\begin{proof}
Put $V=G_{r+1}+tG_{r+1}$. Equation~\eqref{eq:qthree}
implies $xV\subseteq V$, while~\eqref{eq:qvr} gives
\[
 v_r=[r+1]_\zeta^{-1}(v_{r+1}x-\zeta^{r+1}tv_{r+1})\in V,
 \qquad v_rt=v_{r+1}.
\]
Therefore $v_rV\subseteq V$.
The left stabilizer of $V$ is a division subring and hence contains $G_r$, whence $G_r+G_rt\subseteq V$.

Conversely, let $U=G_r+G_rt$.
This is a finite-dimensional left $G_r$-space with $Ux\subseteq U$.
The element $v_{r+1}=v_rt$ belongs to $U$, and
\[
 tv_{r+1}=\zeta^{-r-1}(v_{r+1}x-[r+1]_\zeta v_r)\in U.
\]
Consequently, $Uv_{r+1}\subseteq U$.
The right stabilizer of $U$ is a division subring and contains $G_{r+1}$,
whence $V\subseteq U$. This proves~\eqref{eq:qadjacent}.
Repeated application of~\eqref{eq:qadjacent} gives
\[
 \sum_{i=0}^hG_rt^i
 =\sum_{i=0}^{h-1}(G_r+G_rt)t^i
 =\sum_{i=0}^{h-1}(G_{r+1}+tG_{r+1})t^i
 =\sum_{i=0}^h t^iG_{r+h}.
\]
The last equality follows by induction on $h$.
By Lemma~\ref{lem:qindependence}, the sum on the right-hand side is direct.
\end{proof}

\begin{lemma}\label{lem:qleftbasis}
For $1\le r<N$, we have the decomposition
\begin{equation}\label{eq:qleftbasis}
 E =\bigoplus_{i=0}^{N-r}G_rt^i.
\end{equation}
\end{lemma}

\begin{proof}
Put $\ell=N-r$. By~\eqref{eq:qflag},
\[
 \sum_{i=0}^{\ell-1}G_rt^i
 =\bigoplus_{i=0}^{\ell-1}t^iG_{N-1}.
\]
Lemma~\ref{lem:qindependence} shows that $t^\ell$ does not belong
to this space. If there were a nontrivial left linear relation
among $1,t,\ldots,t^\ell$ with highest exponent $h$,
then normalizing its leading coefficient by multiplication on the left and multiplying on the right
by $t^{\ell-h}$ would give a contradiction.

Let $U=\sum_{i=0}^{\ell}G_rt^i$.
By~\eqref{eq:qpower}, the element $t^N$ belongs to $G_r$, so
\[
 s=v_r(t^N)^{-1}t^\ell\in U.
\]
The formulas $t^ix=\zeta^ixt^i+[i]_\zeta t^{i-1}$ give $Ux\subseteq U$.
Moreover, $t^is=sx^i\in U$, whence $Us\subseteq U$.
The right stabilizer of the finite-dimensional left $G_r$-space $U$
is a division subring and contains $k,x,s$. Consequently,
it equals $E $, and $U=E $.
\end{proof}

\section{The dimension sequence}
\label{sec:dimensions}

\begin{proof}[Proof of Theorem~\ref{thm:quantum}]
Put $H=G_1$. By construction,
\begin{equation}\label{eq:qtwogens}
 E=\Div(x,s),\qquad
 H=\Div(x,st)=\Div(x,s^2xs^{-1}).
\end{equation}
Lemmas~\ref{lem:qrightbasis} and~\ref{lem:qleftbasis} give
\[
 E=H\oplus tH
       =\bigoplus_{i=0}^{N-1}H t^i.
\]
Thus both one-sided dimensions have the required values.

For $0\le j\le N-1$, put
\[
 L_j=\sum_{i=0}^jH t^i.
\]
For $1\le j\le N-2$, define
\[
 S_j=\{a\in E:L_ja\subseteq L_j\},\qquad
 D_j=\{a\in E:L_{j-1}a\subseteq L_j\},
\]
and put $S_0=H$, $D_0=E$.
By~\eqref{eq:qflag} and Lemma~\ref{lem:qrightbasis},
for $0\le j\le N-2$ we have
\[
 L_j=\bigoplus_{i=0}^j t^iG_{j+1},\qquad
 E=\bigoplus_{i=0}^{j+1}t^iG_{j+1}.
\]

For $j\ge1$, the set $S_j$ is a division subring: the right multiplication by any of its nonzero elements is injective,
and hence surjective, on the finite-dimensional left
$H$-space $L_j$.
Moreover, $G_{j+1}\subseteq S_j$.
By the dimension formula, the number
$\dim((S_j)_{G_{j+1}})$ divides both $j+1$ and $j+2$,
so $S_j=G_{j+1}$.

Let $j\ge1$ and $0\ne a\in D_j$.
Since $1\in L_{j-1}$, we have $a\in L_j$.
Write $a=\sum_{i=0}^ht^ib_i$, where
$b_i\in G_{j+1}$, $b_h\ne0$, and $h\le j$.
If $h\ge2$, then $t^{j+1-h}\in L_{j-1}$, but
$t^{j+1-h}a\notin L_j$ by the right linear independence of
$1,t,\ldots,t^{j+1}$ over $G_{j+1}$. Hence $h\le1$.
The reverse inclusion follows from $L_{j-1}t\subseteq L_j$
and the fact that $L_j$ is a right $G_{j+1}$-space.
For $j=0$, we use the decomposition
$E=H\oplus tH$.
Thus
\[
 S_j=G_{j+1},\qquad
 D_j=G_{j+1}\oplus tG_{j+1}
 \qquad(0\le j\le N-2).
\]

Now $1,t,\ldots,t^{N-1}$ form a left $H$-basis of
$E$, and $\dim((D_j)_{S_j})=2$ for $0\le j<N-1$.
By~\cite[Lemma~1.4]{Sar}, the right dimension sequence
has period~\eqref{eq:qperiod}.
It follows from~\cite[Lemma~1.3]{Sar} that $\Lambda_n$ has finite
representation type and exactly $N+2=n$ isomorphism classes of
indecomposable finitely generated right modules.
In the definition of the diagram in~\cite[\S\S\,3--4]{DRS},
the edge label is the length of the period of the corresponding dimension sequence. For~\eqref{eq:qperiod}, this length is $n$, so the diagram is $I_2(n)$.

The diagonal rings of $\Lambda_n$ are division rings, and the off-diagonal
bimodule is projective on both sides; hence $\Lambda_n$
is hereditary. It is left and right artinian, since
${}_H E$ and $E_E$ are finite-dimensional.
\end{proof}

\begin{remark}
The number $N=n-2$ need not be prime: we use the equality
$[N]_\zeta=0$ and the inequalities $[r]_\zeta\ne0$ for $1\le r<N$.
For example, primitive roots of unity of orders $3,4,5,9$ give the diagrams
$I_2(5)$, $I_2(6)$, $I_2(7)$, and $I_2(11)$, respectively.
Both diagonal division rings are infinite-dimensional over $k$, since they contain an element transcendental over $k$, namely $x$.
\end{remark}

\selectlanguage{english}
\end{document}